# Inclusion-Exclusion-Like Identities

Muneerah Al Nuwairan

*Department of mathematics and statistics*

*College of science, King Faisal university, P.O. Box 400*

*Al Ahsa 31982, Saudi Arabia.*

E-mail: msalnuwairan@kfu.edu.sa

**Abstract:**

An expression $E(X_1, \cdots, X_n)$ built using union, intersection, and complements is called inclusion-exclusion-like if, like the union in the exclusion-inclusion principle, there are constants $c_1, \cdots, c_n$ so that for any sequence of sets $\mathcal{A} = (A_1, \cdots, A_n)$ the cardinality of $E(A_1, \cdots, A_n)$ can be expressed as $|E(A_1, \cdots, A_n)| = \Sigma_{k=1}^n c_k i_{n,k}(\mathcal{A})$ where $i_{n,k}(\mathcal{A}) = \sum_{I \subseteq \bar{n}, |I|=k} |\bigcap_{i \in I} A_i|$ is the sum of the cardinalities of all intersections of $k$ sets in the sequence $\mathcal{A}$. In this paper, we construct from the expression $E$, a set of nonempty subsets of $\{1, \cdots, n\}$ called the characteristic set of $E$, and use this set to give a necessary and sufficient condition for the expression to be inclusion-exclusion-like. Furthermore, we give a method for determining the constants $c_1, \cdots, c_n$ in the expression for the cardinality of $E(A_1, \cdots, A_n)$ when it exists. The content of the paper is illustrated by a simple detailed example given in the introduction.

## 1. Introduction

The inclusion-exclusion principle (IEP) is an important tool in combinatorics. It gives a rule for computing the cardinality of the union of a sequence of sets from the sum of the cardinalities of all intersections of $k$ sets in the sequence. Sometimes IEP is stated with the cardinalities being replaced by probabilities [1]. It has many applications in enumeration [4] and [2], cryptography [8], network reliability [6], approximate reasoning [1], and many other fields. In this paper we generalize the classical IEP, viewing it as a rule for giving the cardinality of a set-valued expression $E(\mathcal{A}) = \cup_{i=1}^n A_i$. After listing the notations used and a definition, the results of the paper are outlined. Then we list, for easy reference, the main theorems of the paper and give a simple detailed example of their use. The rest of the paper is devoted to proving these results. Further examples are given as Lemma 3.13 at the end of the paper. Following standard notations, in the following $\bar{n}$ denotes the set of positive integers not exceeding $n$, $\bar{n} = \{1, \cdots, n\}$, $\binom{m}{k}$ denotes the binomial coefficient $\frac{m!}{k!(m-k)!}$ and finally for a set $A$, $|A|$ denotes the cardinality of $A$.

**Definition 1.1**: Let $\mathcal{A} = (A_1, \cdots, A_n)$ be a sequence. For each $1 \leq k \leq n$, we define $i_{n,k}(\mathcal{A})$ to be the sum of the cardinalities of all intersections of $k$ sets in the sequence $\mathcal{A}$. i.e.

$$i_{n,k}(\mathcal{A}) = \sum_{I \subseteq \overline{n}, |I|=k} |\bigcap_{i \in I} A_i|$$

The IEP states that for any sequence $\mathcal{A} = (A_1, \cdots, A_n)$ we have $|A_1 \cup \cdots \cup A_n| = \sum_{k=1}^{n}(-1)^{k-1} i_{n,k}(\mathcal{A})$. In this paper, we determine for an arbitrary expression $E(X_1, \cdots, X_n)$ whether there are constants $c_1, \cdots, c_n$ so that for any sequence $\mathcal{A}$ the cardinality of $E(A_1, \cdots, A_n)$ can be expressed as $|E(A_1, \cdots, A_n)| = \sum_{k=1}^{n} c_k i_{n,k}(\mathcal{A})$. Such expression we call inclusion-exclusion-like (Definition 3.1). In the following, we give a simple algorithm (Theorem 2.11) for computing a unique (Theorem 3.4) set $\mathcal{S} \subseteq \mathcal{P}(\overline{n}) \setminus \{\emptyset\}$ from the expression $E$, and give a simple criterion on $\mathcal{S}$ (Theorem 3.6) that is necessary and sufficient (Theorem 3.11) for the expression $E$ to be inclusion-exclusion-like. In the case where it is, we give a formula for computing the coefficients $c_1, \cdots, c_n$ in the expression of the cardinality (Theorem 3.11). Here is the list of the main results in the paper.

**Theorem 2.11:** For each set-valued expression $E$ there exits a set $\mathcal{S} \subseteq \mathcal{P}(\overline{n}) \setminus \{\emptyset\}$ such that $E \cong F_{\mathcal{S}}$. Furthermore, the set $\mathcal{S}$ can be computed recursively (i.e. as in the definition of set-valued expression, Definition 3.1) as follows

1. If $E = \emptyset$ then $\mathcal{S} = \emptyset$.
2. If $\emptyset \neq I \subseteq \overline{n}$ and $E = \bigcap_{i \in I} X_i$ then $\mathcal{S} = \mathcal{S}_I = \{J \subseteq \overline{n} : I \subseteq J\}$. In particular for $E = X_i$, with $i \in \overline{n}$ we have $\mathcal{S} = \mathcal{S}_{\{i\}}$.
3. If $E = E_1 \cup E_2$ (respectively, $E = E_1 \cap E_2$) and $E_i \cong F_{\mathcal{S}_i}, i = 1,2$ then $E \cong F_{\mathcal{S}}$ where $\mathcal{S} = \mathcal{S}_1 \cup \mathcal{S}_2$ (respectively $\mathcal{S} = \mathcal{S}_1 \cap \mathcal{S}_2$).
4. If $E = E_1^c$ and $E_1 \cong F_{\mathcal{S}_1}$ then $E \cong F_{\mathcal{S}}$ where $\mathcal{S} = \mathcal{P}(\overline{n}) \setminus (\mathcal{S}_1 \cup \{\emptyset\})$.

**Theorem 3.4:** If $\mathcal{S}_1, \mathcal{S}_2 \subseteq \mathcal{P}(\overline{n}) \setminus \{\emptyset\}$ and $F_{\mathcal{S}_1} \cong F_{\mathcal{S}_2}$ then $\mathcal{S}_1 = \mathcal{S}_2$.

**Theorem 3.6:** If the set-valued expression $E$ is inclusion-exclusion-like and $\mathcal{S}$ is its characteristic set, then $\mathcal{S}$ has the following property

> For each $1 \leq k \leq n$, if $\mathcal{S}$ contains a set of cardinality $k$, then it contains every set in $\mathcal{P}(\overline{n})$ of cardinality $k$.

**Theorem 3.11:** Suppose that $E$ is an inclusion-exclusion-like set-valued expression, $\mathcal{S}$ is its characteristic set, and $C = \{k : \exists J \in \mathcal{S} \ni k = |J|\}$, the cardinalities of sets in $\mathcal{S}$. If $f(i) = \chi_C(i)$ is the characteristic function of $C$, then for any $\mathcal{A} = (A_1, \cdots, A_n)$ we have $|E(A_1, \cdots, A_n)| = \sum_{k=1}^{n} c_k i_{n,k}(\mathcal{A})$ where $c_k = \sum_{j=1}^{k}(-1)^{k-j}\binom{k}{j}f(j)$.

An application for using this results is the following:

**Example 1.2:** With $n = 3$, consider two expressions $E = X_1 \cup X_2$ and $F = (X_1 \cap X_2) \cup (X_1 \cap X_3) \cup (X_2 \cap X_3)$.

1. The characteristic set of $E$ is $\mathcal{S}_E = \mathcal{S}_1 \cup \mathcal{S}_2$ (Theorem 2.11 (3)) where $\mathcal{S}_1 = \{\{1\}, \{1,2\}, \{1,3\}, \{1,2,3\}\}$ and $\mathcal{S}_2 = \{\{2\}, \{1,2\}, \{2,3\}, \{1,2,3\}\}$ are the characteristic sets of $X_1$ and $X_2$, respectively (Theorem 2.11 (2))

Since $S_E = \{\{1\}, \{2\}, \{1,2\}, \{1,3\}, \{2,3\}, \{1,2,3\}\}$ contain some subsets of $\overline{3}$ of cardinality 1, namely $\{1\}$ and $\{2\}$ but does not contain all such sets as $\{3\} \notin S_E$ this expression is not inclusion-exclusion-like (Theorem 3.6). That is, one cannot find $c_1, c_2, c_3$ constants so that for any sequence of sets $\mathcal{A} = (A_1, A_2, A_3)$ we have that $|A_1 \cup A_2| = c_1 i_{3,1}(\mathcal{A}) + c_2 i_{3,2}(\mathcal{A}) + c_3 i_{3,3}(\mathcal{A})$. This can be seen directly by considering, as in the proof of Theorem 3.6, the two sequences (Definition 3.2) $^{\{2\}}\mathcal{A} = (\emptyset, \{1\}, \emptyset)$ and $^{\{3\}}\mathcal{A} = (\emptyset, \emptyset, \{1\})$. We clearly have
$$i_{3,k}(^{\{2\}}\mathcal{A}) = i_{3,k}(^{\{3\}}\mathcal{A}) = \begin{cases} 1 & k = 1 \\ 0 & k \neq 1 \end{cases}$$
so any expression of the type $|A_1 \cup A_2| = c_1 i_{3,1}(\mathcal{A}) + c_2 i_{3,2}(\mathcal{A}) + c_3 i_{3,3}(\mathcal{A})$ will give us that $1 = |^{\{2\}}A_1 \cup {}^{\{2\}}A_2| = \Sigma_{k=1}^3 c_k i_{3,k}(^{\{2\}}\mathcal{A}) = \Sigma_{k=1}^3 c_k i_{3,k}(^{\{3\}}\mathcal{A}) = |^{\{3\}}A_1 \cup {}^{\{3\}}A_2| = 0$

2. As in 1), using Theorem 2.11 (2) and (3), we find that the characteristic set for $F$ is $S_F = \{\{1,2\}, \{1,3\}, \{2,3\}, \{1,2,3\}\}$. This consists of all subsets of $\overline{3}$ with cardinality falling in the set $C = \{2,3\}$ so (by Theorem 3.11) there are constants $c_1, c_2, c_3$ so that $|F(A_1 A_2, A_3)| = c_1 i_{3,1}(\mathcal{A}) + c_2 i_{3,2}(\mathcal{A}) + c_3 i_{3,3}(\mathcal{A})$ for any sequence of sets $\mathcal{A} = (A_1, A_2, A_3)$. The constants $c_1, c_2, c_3$ can be computed as in Theorem 3.11 by taking $f(n) = \chi_{\{2,3\}}$, the characteristic function of $C$, and $c_k = \Sigma_{j=1}^k (-1)^{k-j} \binom{k}{j} f(j)$.
That gives $c_1 = 0, c_2 = 1, c_3 = -2$. We obtain the identity $|(A_1 \cap A_2) \cup (A_1 \cap A_3) \cup (A_2 \cap A_3)| = i_{3,2}(\mathcal{A}) - 2\, i_{3,3}(\mathcal{A})$ which is a special case of the generalized inclusion-exclusion principle (Lemma 3.13 1)).

## 2. Characteristic sets

In the following, we keep $n \in N$ constant and use the variables $X_1, \cdots, X_n$. In this section, we construct for each equivalence class of set-valued expression $E$ a set $S_E \subseteq \mathcal{P}(\overline{n})$. This construction is analogous to the disjunctive normal form in Boolean algebras [10, section 2.4]. We start with some definitions.

**Definition 2.1:** Set-valued expressions in the variables $X_1, \cdots, X_n$ are defined recursively as follows:

1. $\emptyset$ is a (constant) set-valued function.
2. For each $1 \leq i \leq n$, $X_i$ is a set-valued function.
3. If $F(X_1, \cdots, X_n)$ is a set-valued expression, then so is $F(X_1, \cdots, X_n)^c$.
4. If $F_1(X_1, \cdots, X_n)$ and $F_2(X_1, \cdots, X_n)$ are set-valued expressions, then so are both $F_1(X_1, \cdots, X_n) \cup F_1(X_1, \cdots, X_n)$ and $F_1(X_1, \cdots, X_n) \cap F_1(X_1, \cdots, X_n)$.

Given a sequence $\mathcal{A} = (A_1, \cdots, A_n)$ of sets, we define the value of a set-valued expression to be the set resulting from substituting $A_i$ for $X_i$ and evaluating the result with complements being taken with respect to $A_1 \cup \cdots \cup A_n$. This is given explicitly in the following definition:

**Definition 2.2:** Given a sequence of sets $\mathcal{A} = (A_1, \cdots, A_n)$ the value of a set-valued expression $F(X_1, \cdots, X_n)$ at $\mathcal{A}$ denoted $F(\mathcal{A})$ or $F(A_1, \cdots, A_n)$ is defined recursively as follows:

1. If $F = \emptyset$ then $F(\mathcal{A}) = \emptyset$.

2. If $F = X_i, 1 \leq i \leq n$ then $F(\mathcal{A}) = A_i$.
3. If $F = E(X_1, \cdots, X_n)^c$ then $F(\mathcal{A}) = (A_1 \cup \cdots \cup A_n) \setminus E(\mathcal{A})$.
4. If $F = E_1 \cup E_2$ (resp. $F = E_1 \cap E_2$) then $F(\mathcal{A}) = E_1(\mathcal{A}) \cup E_2(\mathcal{A})$ (resp $F(\mathcal{A}) = E_1(\mathcal{A}) \cap E_2(\mathcal{A})$.

**Definition 2.3:** Two set-valued expressions $F(X_1, \cdots, X_n)$, and $E(X_1, \cdots, X_n)$ are equivalent, denoted $F(X_1, \cdots, X_n) \cong E(X_1, \cdots, X_n)$, if for all sequences of sets $\mathcal{A} = (A_1, \cdots, A_n)$ we have $F(\mathcal{A}) = E(\mathcal{A})$.

The relation $\cong$ is easily seen to be an equivalence relation satisfying the following properties

**Lemma 2.4:** For set valued expressions $E_1, E_2, E_3, F_1, F_2$ we have:

1. If $E_1 \cong F_1$ and $E_2 \cong F_2$ then $E_1 \cup E_2 \cong F_1 \cup F_2$ and $E_1 \cap E_2 \cong F_1 \cap F_2$ and $E_1^c \cong F_1^c$.
2. $E_1 \cup E_2 \cong E_2 \cup E_1$ and $E_1 \cap E_2 \cong E_2 \cap E_1$.
3. $E_1 \cup (E_2 \cup E_3) \cong (E_1 \cup E_2) \cup E_3$ and $E_1 \cap (E_2 \cap E_3) \cong (E_1 \cap E_2) \cap E_3$.
4. $E_1 \cup E_1 \cong E_1$, $E_1 \cap E_1 \cong E_1$, $E_1 \cup \emptyset \cong E_1$, and $E_1 \cap \emptyset \cong \emptyset$.

**Proof:**

Let $\mathcal{A}$ be a sequence of sets. The premise gives $E_1(\mathcal{A}) = F_1(\mathcal{A})$ and $E_2(\mathcal{A}) = F_2(\mathcal{A})$. Thus $(E_1 \cup E_2)(\mathcal{A}) = E_1(\mathcal{A}) \cup E_2(\mathcal{A}) = F_1(\mathcal{A}) \cup F_2(\mathcal{A}) = (F_1 \cup F_2)(\mathcal{A})$. Since $\mathcal{A}$ was arbitrary $E_1 \cup E_2 \cong F_1 \cup F_2$. The proof for $E_1 \cap E_2 \cong F_1 \cap F_2$ and $E_1^c \cong F_1^c$. The rest are similar. ∎

**Remark 2.5:**

Lemma 2.4 (1.) gives us that we can define union, intersection and complement unambiguously on the equivalence classes. 2. and 3. show that union and intersection are commutative and associative operations on the equivalence classes. Thus if $E_i, i \in I$ are expressions then $\bigcap_{i \in I} E_i$ is well-defined up to an equivalence. That is, the way in which we take the intersection of $E_i, i \in I$ does not alter the equivalence class of the result. Likewise, $\bigcup_{i \in I} E_i$ is well-defined up to an equivalence and $E^c$ for an expression $E$. In the following we will, without explicit mention, use arbitrary unions and intersections with the understanding that they are defined up to an equivalence.

**Definition 2.6:**

1. For $I \subseteq \overline{n}$ we define $G_I(X_1, \cdots, X_n) = \bigcap_{i \in I} X_i \cap \bigcap_{i \in \overline{n} \setminus I} X_i^c$.
2. For $\mathcal{S} \subseteq \mathcal{P}(\overline{n})$ we define $F_{\mathcal{S}}(X_1, \cdots, X_n) = \bigcup_{I \in \mathcal{S}} G_I(X_1, \cdots, X_n)$.

We start by listing some useful properties of the functions defined above

**Proposition 2.7:** The set-valued functions $G_I, I \subseteq \overline{n}$ satisfy the following

1. For any sequence $\mathcal{A} = (A_1, \cdots, A_n)$ and $x \in \bigcup_{i \in \overline{n}} A_i$, we have $x \in G_I(\mathcal{A})$ if and only if $\{k: x \in A_k\} = I$.
2. $G_\emptyset \cong \emptyset$.
3. If $I_1 \neq I_2$ then $G_{I_1}(X_1, \cdots, X_n) \cap G_{I_2}(X_1, \cdots, X_n) \cong \emptyset$.
4. $\bigcup_{I \in \mathcal{P}(\overline{n})} G_I(X_1, \cdots, X_n) \cong X_1 \cup \cdots \cup X_n$.

**Proof:**

1. $x \in G_I(\mathcal{A})$ if and only if $x \in \bigcap_{i \in I} A_i \cap \bigcap_{i \in \bar{n} \setminus I} A_i^c$ if and only if $x \in A_i$ for each $i \in I$ and $x \notin A_i$ for each $i \in \bar{n} \setminus I$. That is if and only $\{k: x \in A_k\} = I$
2. $G_\emptyset(X_1, \cdots, X_n) = \bigcap_{i \in \bar{n}} X_i^c$ so for any sequence $\mathcal{A}$ we have by definition of value and DeMorgan law $G_\emptyset(\mathcal{A}) = \bigcap_{i \in \bar{n}} A_i^c = (\bigcup_{i \in \bar{n}} A_i)^c = \bigcup_{i \in \bar{n}} A_i \setminus \bigcup_{i \in \bar{n}} A_i = \emptyset$. The result follows by definition of equivalence
3. For any sequence $\mathcal{A} = (A_1, \cdots, A_n)$, if $x \in G_{I_1}(A_1, \cdots, A_n)$ then by 1) we have $\{k: x \in A_k\} = I_1 \neq I_2$ so by 1) again $x \notin G_{I_2}(A_1, \cdots, A_n)$. Thus $G_{I_1}(\mathcal{A}) \cap G_{I_2}(\mathcal{A}) = \emptyset$. Since $\mathcal{A}$ was an arbitrary sequence we have $G_{I_1}(X_1, \cdots, X_n) \cap G_{I_2}(X_1, \cdots, X_n) \cong \emptyset$
4. Let $\mathcal{A} = (A_1, \cdots, A_n)$ be a sequence of sets. It is clear (e.g. by induction) that for any set valued expression $E(X_1, \cdots, X_n)$ we have $E(\mathcal{A}) \subseteq A_1 \cup \cdots \cup A_n$ so $(\bigcup_{I \in \mathcal{P}(\bar{n})} G_I)(\mathcal{A}) \subseteq A_1 \cup \cdots \cup A_n$. On the other hand, if $x \in A_1 \cup \cdots \cup A_n$ and we take $I = \{k: x \in A_k\}$ then $I \neq \emptyset$ and, by (1), $x \in G_I(\mathcal{A})$ so $x \in (\bigcup_{I \in \mathcal{P}(\bar{n})} G_I)(\mathcal{A})$. Thus $(\bigcup_{I \in \mathcal{P}(\bar{n})} G_I)(\mathcal{A}) = A_1 \cup \cdots \cup A_n$. Since $\mathcal{A}$ was arbitrary we have $\bigcup_{I \in \mathcal{P}(\bar{n})} G_I(X_1, \cdots, X_n) \cong X_1 \cup \cdots \cup X_n$. ∎

**Definition 2.8:** For $I \subseteq \bar{n}$ we define $\mathcal{S}_I \subseteq \mathcal{P}(\bar{n})$ by $\mathcal{S}_I = \{J \subseteq \bar{n}: I \subseteq J\}$.

**Proposition 2.9:** The set-valued expression $F_\mathcal{S}(X_1, \cdots, X_n)$ satisfy the following properties:

1. $F_\emptyset \cong F_{\{\emptyset\}} \cong \emptyset$ and $F_{\mathcal{P}(\bar{n})} \cong X_1 \cup \cdots \cup X_n$.
2. For any $\mathcal{S} \subseteq \mathcal{P}(\bar{n})$ we have $F_\mathcal{S} \cong F_{\mathcal{S} \cup \{\emptyset\}}$.
3. For $\emptyset \neq I \subseteq \bar{n}$ we have $F_{\mathcal{S}_I}(X_1, \cdots, X_n) \cong \bigcap_{i \in I} X_i$.
4. $F_{\mathcal{S}_1} \cup F_{\mathcal{S}_2} \cong F_{\mathcal{S}_1 \cup \mathcal{S}_2}$.
5. $F_{\mathcal{S}_1} \cap F_{\mathcal{S}_2} \cong F_{\mathcal{S}_1 \cap \mathcal{S}_2}$.
6. $F_\mathcal{S}^c \cong F_{\mathcal{P}(\bar{n}) \setminus (\mathcal{S} \cup \{\emptyset\})}$.

**Proof:**

1. For any $\mathcal{A} = (A_1, \cdots, A_n)$ we have $F_\emptyset(\mathcal{A}) = \emptyset$, since it is an empty union. Also, by definition of $F_{\{\emptyset\}}$ and Proposition 2.7(2), we have $F_{\{\emptyset\}}(X_1, \cdots, X_n) = G_\emptyset(X_1, \cdots, X_n) \cong \emptyset$. The last equivalence follows from Proposition 2.7 (4).
2. If $\emptyset \in \mathcal{S}$ the two sides are identical ($\mathcal{S} = \mathcal{S} \cup \{\emptyset\}$). If $\emptyset \notin \mathcal{S}$ then $F_{\mathcal{S} \cup \{\emptyset\}} = F_\mathcal{S} \cup G_\emptyset$ (by definition of $F_{\mathcal{S} \cup \{\emptyset\}}$). But $G_\emptyset \cong \emptyset$ so $F_{\mathcal{S} \cup \{\emptyset\}} \cong F_\mathcal{S} \cup \emptyset \cong F_\mathcal{S}$.
3. Let $\mathcal{A} = (A_1, \cdots, A_n)$ be a sequence of sets. Since $I \subseteq J$ for each $J \in \mathcal{S}_I$ we have $\bigcap_{i \in I} A_i \supseteq \bigcap_{i \in J} A_i \supseteq \bigcap_{i \in J} A_i \cap \bigcap_{i \in \bar{n} \setminus J} A_i^c = G_J(\mathcal{A})$. Thus $\bigcap_{i \in I} A_i \supseteq (\bigcup_{J \in \mathcal{S}_I} G_J)(\mathcal{A})$ On the other hand, if $x \in \bigcap_{i \in I} A_i$ and $J = \{k: x \in A_k\}$ then $J \supseteq I$, i.e. $J \in \mathcal{S}_I$ and by Proposition 2.7(1) $x \in G_J(\mathcal{A}) \subseteq (\bigcup_{J \in \mathcal{S}_I} G_J)(\mathcal{A})$. Thus $\bigcap_{i \in I} A_i = (\bigcup_{J \in \mathcal{S}_I} G_J)(\mathcal{A}) = F_{\mathcal{S}_I}(\mathcal{A})$. Since $\mathcal{A}$ is arbitrary this gives us that $F_{\mathcal{S}_I}(X_1, \cdots, X_n) \cong \bigcap_{i \in I} X_i$.
4. Follows easily from definition and Lemma 2.4 (4).
5. Let $\mathcal{A} = (A_1, \cdots, A_n)$ be a sequence of sets. If $x \in F_{\mathcal{S}_1} \cap F_{\mathcal{S}_2}(\mathcal{A})$ then there exist sets $J_1 \in \mathcal{S}_1$ and $J_2 \in \mathcal{S}_2$ such that $x \in G_{J_1}(\mathcal{A})$ and $x \in G_{J_2}(\mathcal{A})$. But then $G_{J_1}(\mathcal{A}) \cap G_{J_2}(\mathcal{A}) \neq \emptyset$ so by Proposition 2.7(2) $J_1 = J_2$ and this set is in $\mathcal{S}_1 \cap \mathcal{S}_2$.

Thus $x \in G_{J_1}(\mathcal{A}) \subseteq F_{S_1 \cap S_2}(\mathcal{A})$ and $F_{S_1} \cap F_{S_2}(\mathcal{A}) \subseteq F_{S_1 \cap S_2}(\mathcal{A})$. The opposite containment follows since $F_{S_1 \cap S_2}(\mathcal{A}) = \cup_{I \in S_1 \cap S_2} G_I(\mathcal{A}) \subseteq \left( \cup_{I \in S_1} G_I(\mathcal{A}) \right) \cap \left( \cup_{I \in S_2} G_I(\mathcal{A}) \right) = F_{S_1} \cap F_{S_2}(\mathcal{A})$ giving us $F_{S_1} \cap F_{S_2}(\mathcal{A}) = F_{S_1 \cap S_2}(\mathcal{A})$. Since $\mathcal{A}$ is arbitrary $F_{S_1} \cap F_{S_2} \cong F_{S_1 \cap S_2}$.

6. For $\mathcal{A} = (A_1, \cdots, A_n)$ by 5. and 1. $F_S(\mathcal{A}) \cap F_{\mathcal{P}(\overline{n}) \setminus (S \cup \{\emptyset\})}(\mathcal{A}) = F_\emptyset(\mathcal{A}) = \emptyset$. We also have from 3. that $F_{\mathcal{P}(\overline{n}) \setminus (S \cup \{\emptyset\})}(\mathcal{A}) \cup F_S(\mathcal{A}) = F_{(\mathcal{P}(\overline{n}) \setminus (S \cup \{\emptyset\})) \cup S}(\mathcal{A})$. This is $F_{\mathcal{P}(\overline{n})}(\mathcal{A})$ (when $\emptyset \in S$) or $F_{\mathcal{P}(\overline{n}) \setminus \{\emptyset\}}(\mathcal{A})$ (when $\emptyset \notin S$). But the two are equal by 2. and by 1. both equal $A_1 \cup \cdots \cup A_n$. Thus $F_{\mathcal{P}(\overline{n}) \setminus (S \cup \{\emptyset\})}(\mathcal{A}) = F_S(\mathcal{A})^c$. ∎

**Remark 2.10:**

Our goal is to establish a bijective correspondence between equivalence classes of set-valued functions and subsets of $\mathcal{P}(\overline{n})$ where if $S$ corresponds to the equivalence class of $E$ then $E \cong F_S$. Without further restriction, the class $S$ would not be unique since $F_\emptyset \cong F_{\{\emptyset\}}$. It turns out this non-uniqueness can be avoided by taking $S$ to be a (possibly empty) collection of nonempty subsets of $\overline{n}$.

**Theorem 2.11:** For each set-valued expression $E$ there exits a set $S \subseteq \mathcal{P}(\overline{n}) \setminus \{\emptyset\}$ such that $E \cong F_S$. Furthermore, the set $S$ can be computed recursively (i.e. as in the definition of set-valued expression) as follows

1. If $E = \emptyset$ then $S = \emptyset$.
2. If $\emptyset \neq I \subseteq \overline{n}$ and $E = \cap_{i \in I} X_i$ then $S = S_I = \{J \subseteq \overline{n}: I \subseteq J\}$. In particular for $E = X_i$, with $i \in \overline{n}$ we have $S = S_{\{i\}}$.
3. If $E = E_1 \cup E_2$ (respectively $E = E_1 \cap E_2$) and $E_i \cong F_{S_i}, i = 1,2$ then $E \cong F_S$ where $S = S_1 \cup S_2$ (respectively $S = S_1 \cap S_2$).
4. If $E = E_1^c$ and $E_1 \cong F_{S_1}$ then $E \cong F_S$ where $S = \mathcal{P}(\overline{n}) \setminus (S_1 \cup \{\emptyset\})$.

**Proof:**

From Proposition 2.9 (1 and 3) we see that if $E = \emptyset$ (respectively $E = \cap_{i \in I} X_i$ where $\emptyset \neq I \subseteq \overline{n}$) then $E \cong F_S$ where $S$ is $\emptyset$ (respectively $S_I$) establishing 1. and 2. Also, items 3 and 4 follow from Proposition 2.9 (3 and 4). We finally note that in the cases 1, 2, and 4, the set $S$ does not contain the empty set and the other cases are built from previous stages by using union or intersection so will not contain the empty set if the previous stages don't. ∎

### 3. Inclusion-Exclusion-like expressions

**Definition 3.1:** A set-valued expression is called inclusion-exclusion-like if there are constants $c_1, \cdots, c_n$ such that for any sequence of sets $\mathcal{A} = (A_1, \cdots, A_n)$ we have $|E(\mathcal{A})| = \Sigma_{k=1}^n c_k i_{n,k}(\mathcal{A})$.

The sequences in the following definition will be useful in our study of inclusion-exclusion-like expressions.

**Definition 3.2:** For $\emptyset \neq I \subseteq \overline{n}$, we take $^I\mathcal{A} = (^IA_1, \cdots, ^IA_n)$ where

$$^IA_k = \begin{cases} \{1\} & k \in I \\ \emptyset & k \notin I \end{cases}$$

We start with listing some properties of the sequences $^I\mathcal{A}$.

**Proposition 3.3:** For each $\emptyset \neq I \subseteq \overline{n}$ we have:

1. $\cap_{i \in J} {}^I A_i = \begin{cases} \{1\} & if \ J \subseteq I \\ \emptyset & if \ J \not\subseteq I \end{cases}$
2. $i_{n,k}(^I\mathcal{A}) = \binom{|I|}{k}$.
3. $G_J(^I\mathcal{A}) = \{1\}$ if $J = I$ and $G_J(^I\mathcal{A}) = \emptyset$ if $J \neq I$.
4. For $\mathcal{S} \subseteq \mathcal{P}(\overline{n}) \setminus \{\emptyset\}$ we have $F_\mathcal{S}(^I\mathcal{A}) = \begin{cases} \{1\} & if \ I \in \mathcal{S} \\ \emptyset & if \ I \notin \mathcal{S} \end{cases}$

**Proof:**

Item 1. is clear. From 1, we have $i_{n,k}(^I\mathcal{A}) = \Sigma_{|J|=k} |\cap_{i \in J} {}^IA_i| = \Sigma_{|J|=k, J \subseteq I} 1 = \binom{|I|}{k}$ (the number of subsets of $I$ of size $k$) which proves 2. For 3, note that for $^I\mathcal{A}$ the complement is taken with respect to $\cup_{i \in \overline{n}} {}^IA_i = \{1\}$, so we have

$$^IA_k^c = \{1\} \setminus {}^IA_k = \begin{cases} \emptyset & if \ k \in I \\ \{1\} & if \ k \in \overline{n} \setminus I \end{cases}.$$

Thus, we have $G_I(^I\mathcal{A}) = \cap_{i \in I} {}^IA_i \cap \cap_{i \in \overline{n} \setminus I} {}^IA_i^c = \{1\}$. If $J \neq I$ then either there exists $j \in J$ with $j \in \overline{n} \setminus I$ or there exist $i \in I$ with $i \in \overline{n} \setminus J$. In the first case we have $G_J(^I\mathcal{A}) \subseteq \cap_{i \in J} {}^IA_i \subseteq {}^IA_j = \emptyset$, in the second $G_J(^I\mathcal{A}) \subseteq \cap_{j \in \overline{n} \setminus J} {}^IA_j^c \subseteq {}^IA_i^c = \emptyset$. In both cases we have $G_J(^I\mathcal{A}) = \emptyset$. Item 4, follows immediately from 3 and the definition of $F_\mathcal{S}$. ∎

Using Proposition 3.3, we obtain the uniqueness of the set $\mathcal{S}$ such that $E \cong F_\mathcal{S}$.

**Theorem 3.4:** If $\mathcal{S}_1, \mathcal{S}_2 \subseteq \mathcal{P}(\overline{n}) \setminus \{\emptyset\}$ and $F_{\mathcal{S}_1} \cong F_{\mathcal{S}_2}$ then $\mathcal{S}_1 = \mathcal{S}_2$.

**Proof:**

Suppose $\mathcal{S}_1 \neq \mathcal{S}_2$ without loss of generality there is $I \in \mathcal{S}_1$ such that $I \notin \mathcal{S}_2$. From proposition 3.3 we have $F_{\mathcal{S}_1}(^I\mathcal{A}) = \{1\} \neq \emptyset = F_{\mathcal{S}_2}(^I\mathcal{A})$ so $F_{\mathcal{S}_1} \not\cong F_{\mathcal{S}_2}$. ∎

**Definition 3.5:** For a set-valued expression $E$ the unique set $\mathcal{S} \subseteq \mathcal{P}(\overline{n}) \setminus \{\emptyset\}$ such that $E \cong F_\mathcal{S}$ is called the characteristic set of $E$.

Now we deduce a necessary and sufficient condition for a set-valued expression to be inclusion-exclusion like.

**Theorem 3.6:** If the set-valued expression $E$ is inclusion-exclusion-like and $\mathcal{S}$ is its characteristic set then $\mathcal{S}$ has the following property:

> For each $1 \leq k \leq n$, if $\mathcal{S}$ contains a set of cardinality $k$ then it contains every set in $\mathcal{P}(\overline{n})$ of cardinality $k$.

**Proof:**

The theorem's assumption is equivalent to assuming that $F_S$ is inclusion-exclusion-like. That is, there exists a sequence of numbers $c_1, \cdots, c_n$ such that for any sequence $\mathcal{A} = (A_1, \cdots, A_n)$ of sets we have (∗) $|F_S(\mathcal{A})| = c_1 i_{n,1}(\mathcal{A}) + \cdots + c_n i_{n,n}(\mathcal{A})$.

Now suppose that $I \in S$ has cardinality $k$ and that $J \in \mathcal{P}(\overline{n})$ also has cardinality $k$. By Proposition 3.3 (4) we have $|F_S(^I\mathcal{A})| = 1$ so by (∗) and Proposition 3.3 (2) we have

$1 = c_1 i_{n,1}(^I\mathcal{A}) + \cdots + c_n i_{n,n}(^I\mathcal{A}) = c_1 \binom{k}{1} + c_2 \binom{k}{2} + \cdots + c_n \binom{k}{n}$. Applying Proposition 3.3 (2) and (∗) with $J$ replacing $I$ we obtain $|F_S(^J\mathcal{A})| = c_1 \binom{k}{1} + c_2 \binom{k}{2} + \cdots + c_n \binom{k}{n}$. Thus, from the previous equality $|F_S(^J\mathcal{A})| = 1$. From this and Proposition 3.3 (4) we have $J \in S$. ∎

Our next goal is to show that the above necessary condition for a set-valued expression to be inclusion-exclusion-like is also sufficient. To prove this, we need some combinatorial identities. These identities have been well-known for a long time (see e.g. [5, p.153]). For completeness we add proofs.

**Definition 3.7 :** For $1 \leq k \leq n$

1) $S_k(X_1, \cdots, X_n) = \bigcup_{|I|=k} G_I(X_1, \cdots, X_n)$
2) For a sequence $\mathcal{A} = (A_1, \cdots, A_n)$ we define $\sigma_{n,k}(\mathcal{A}) = |S_k(\mathcal{A})| = \sum_{|I|=k} |G_I(\mathcal{A})|$
   (where the last equality is an immediate consequence of Proposition 2.7 (3))

**Proposition 3.8:** Let $\mathcal{A} = (A_1, \cdots, A_n)$ be a sequence of sets

1. For each nonempty $I \subseteq \mathcal{P}(\overline{n})$ we have $|\bigcap_{i \in I} A_i| = \sum_{I \subseteq J} |G_J(\mathcal{A})|$
2. For $1 \leq k \leq n$ we have $i_{n,k}(\mathcal{A}) = \sum_{j=k}^{n} \binom{j}{k} \sigma_{n,j}(\mathcal{A})$
3. $\sum_{j=1}^{n} \sigma_{n,j}(\mathcal{A}) = |A_1 \cup \cdots \cup A_n|$

**Proof:**

1. By propositions 2.9 (3) and 2.7(3) the definition of $F_S$ for a set $S$

$$\left|\bigcap_{i \in I} A_i\right| = |F_{\{J: I \subseteq J \subseteq \overline{n}\}}(\mathcal{A})| = \left|\bigcup_{J \supseteq I} G_J(\mathcal{A})\right| = \sum_{I \subseteq J} |G_J(\mathcal{A})|$$

2. By item 1, and the definition of $i_{n,k}$ we have

$$i_{n,k}(\mathcal{A}) = \sum_{|I|=k} \left|\bigcap_{i \in I} A_i\right| = \sum_{|I|=k} \sum_{I \subseteq J} |G_J(\mathcal{A})|$$
$$= \sum_{|J| \geq k} \sum_{\{I: I \subseteq J, |I|=k\}} |G_J(\mathcal{A})|$$

In the last sum $|G_J(\mathcal{A})|$ will occur as many times as there are subsets $I$ of $J$ of cardinality $k$, i.e. $\binom{|J|}{k}$. So the last sum is $\sum_{|J| \geq k} \binom{|J|}{k} |G_J(\mathcal{A})|$ or separating this sums into sums based on the cardinality of $J$ we obtain $\sum_{j \geq k} \binom{j}{k} \sum_{|J|=j} |G_J(\mathcal{A})|$ which, by definition of $\sigma_{n,k}(\mathcal{A})$, is $\sum_{j=k}^{n} \binom{j}{k} \sigma_{n,j}(\mathcal{A})$ proving 2.

3. $\sum_{j=1}^n \sigma_{n,j}(\mathcal{A}) = \sum_{j=1}^n \sum_{|I|=j} |G_I(\mathcal{A})| = \sum_{\emptyset \neq I \subseteq \overline{n}} |G_I(\mathcal{A})| = |A_1 \cup \cdots \cup A_n|$, where the last equality follows from Proposition 2.7 (3 and 4). ∎

As we see below, to obtain an expression for the cardinality of $|E(\mathcal{A})|$ in terms of $i_{n,1}, i_{n,2}, \cdots, i_{n,n}$ when $E$ satisfies the conditions in Theorem 3.6, we need an expression for $\sigma_{n,k}(\mathcal{A})$ in terms of $i_{n,1}, i_{n,2}, \cdots, i_{n,n}$. Such an expression can be derived directly from definitions. Here we derive it from Proposition 3.8 using a general principle of inversion.

The following proposition is obtained from [9, p.45 equation (4)] or [7 p.192 formula 5.48 with $f(k) = (-1)^k b_k, g(k) = a_k$ ].

**Proposition 3.9:** If $a_0, \cdots, a_n$ and $b_0, \cdots, b_n$ are sequences of numbers, then

$$b_k = \sum_{j=k}^n (-1)^{j-k} \binom{j}{k} a_j \quad \text{if and only if} \quad a_k = \sum_{j=k}^n \binom{j}{k} b_j.$$

As a corollary to this proposition and Proposition 3.8, we obtain

**Corollary 3.10:** For any sequence $\mathcal{A} = (A_1, \cdots, A_n)$ we have

$$\sigma_{n,k}(\mathcal{A}) = \sum_{j=k}^n (-1)^{j-k} \binom{j}{k} i_{n,j}(\mathcal{A})$$

**Proof:**

Let $b_k = \sigma_{n,k}(\mathcal{A})$ and $a_k = i_{n,k}(\mathcal{A})$ for $1 \leq k \leq n$ and take $a_0 = |A_1 \cup \cdots A_n|$ and $b_0 = 0$. We have, by Proposition 3.8: 2) and 3), $a_k = \sum_{j=k}^n \binom{j}{k} b_j$ for $0 \leq k \leq n$. Thus, by Proposition 3.9, $\sigma_{n,k}(\mathcal{A}) = \sum_{j=k}^n (-1)^{j-k} \binom{j}{k} i_{n,j}(\mathcal{A})$ for $1 \leq k \leq n$. ∎

Using this identity, we show that the condition in Theorem 3.4 is also sufficient.

**Theorem 3.11:** Suppose that $E$ is an inclusion-exclusion-like set-valued expression, $\mathcal{S}$ is its characteristic set, and $C = \{k: \exists J \in \mathcal{S} \ni k = |J|\}$, the cardinalities of sets in $\mathcal{S}$. If $f(i) = \chi_C(i)$ is the characteristic function of $C$, then for any $\mathcal{A} = (A_1, \cdots, A_n)$, we have $|E(A_1, \cdots, A_n)| = \sum_{k=1}^n c_k i_{n,k}(\mathcal{A})$ where $c_k = \sum_{j=1}^k (-1)^{k-j} \binom{k}{j} f(j)$.

**Proof:**

The condition of Theorem 3.6 gives us $\mathcal{S} = \{I \in \mathcal{P}(\overline{n}): |I| \in C\}$, where $C$ is as in the theorem. By definition of characteristic set $E(\mathcal{A}) = F_\mathcal{S}(\mathcal{A}) = \bigcup_{|I| \in C} G_I(\mathcal{A}) = \bigcup_{k \in C} \bigcup_{|I|=k} G_I(\mathcal{A})$. Since the sets $G_I(\mathcal{A}), I \in \mathcal{P}(\overline{n})$ are pairwise disjoint we have

$$|E(\mathcal{A})| = |F_\mathcal{S}(\mathcal{A})| = \sum_{j \in C} \sum_{|I|=j} |G_I(\mathcal{A})| = \sum_{j \in C} \sigma_{n,j}(\mathcal{A}) = \sum_{j=1}^n f(j) \sigma_{n,j}(\mathcal{A})$$

From corollary 3.10 we have

$|E(\mathcal{A})| = \sum_{j=1}^{n} f(j) \sum_{k=j}^{n} (-1)^{k-j} \binom{k}{j} i_{n,k(\mathcal{A})} = \sum_{1 \leq j \leq n, j \leq k \leq n} f(j) \binom{k}{j} (-1)^{k-j} i_{n,k(\mathcal{A})} =$

$\sum_{1 \leq k \leq n, 1 \leq j \leq k} f(j) \binom{k}{j} (-1)^{k-j} i_{n,k}(\mathcal{A}) = \sum_{k=1}^{n} c_k \, i_{n,k(\mathcal{A})}$

Where $c_k = \sum_{j=1}^{k} (-1)^{k-j} \binom{k}{j} f(j)$. ∎

**Remark 3.12:** By proposition 2.7 (1) and the definition of $F_S$, for any sequence $\mathcal{A} = (A_1, \cdots, A_n)$, we have $x \in F_S(\mathcal{A})$ if and only if $\{k: x \in A_k\} \in S$.

In particular, for $S = \{I \subseteq \mathcal{P}(\overline{n}) \setminus \{\emptyset\}: |I| \in C\}$, the characteristic set of an inclusion-exclusion-like expression (where $C \subseteq \overline{n}$), we have $F_S(\mathcal{A}) = \{x: |\{k: x \in A_k\}| \in C\}$ In the lemma below we give three examples of the use of Theorem 3.11 using the, possibly more familiar, right hand side of the last equality instead of $F_S(\mathcal{A})$.

**Lemma 3.13:** Let $n \in \mathbb{N}$. For any sequence $\mathcal{A} = (A_1, \cdots, A_n)$ of sets, we have

1. (The generalized inclusion-exclusion principle)

$$|\{x: |\{j: x \in A_j\}| \geq m\}| = \sum_{k=m}^{n} (-1)^{k-m} \binom{k-1}{k-m} i_{n,k}(\mathcal{A})$$

2. $|\{x: |\{j: x \in A_j\}| \text{ is even}\}| = \sum_{k=1}^{n} (-1)^k (2^{k-1} - 1) i_{n,k}(\mathcal{A})$.
3. $|\{x: |\{j: x \in A_j\}| \text{ is odd}\}| = \sum_{k=1}^{n} (-1)^{k-1} 2^{k-1} i_{n,k}(\mathcal{A})$.

**Proof:**

As noted in the remark, the above expressions correspond to $F_S(\mathcal{A})$ where

$S = \{I \in \mathcal{P}(\overline{n}) \setminus \{\emptyset\}: |I| \in C\}$ and the set $C$ in the three cases is

(1) $\{k: m \leq k \leq n\}$     (2) $\{2k : 1 \leq k \leq n/2\}$     (3) $\{2k + 1: 0 \leq k \leq (n-1)/2\}$

Thus by Theorem 3.12, we have in the three cases $|F_S(\mathcal{A})| = \sum_{k=1}^{n} c_k i_{n,k}(\mathcal{A})$ where $c_k = \sum_{j=1}^{k} (-1)^{k-j} \binom{k}{j} f(j)$ with $f$ the characteristic function of $C$. This gives us:

For 1, $c_k = \sum_{j=m}^{k} (-1)^{k-j} \binom{k}{j}$. Thus, $c_k = 0$ for $k < m$, and for $k \geq m$

$c_k = \sum_{j=m}^{k} (-1)^{k-j} \binom{k}{j} = \sum_{j=m}^{k} (-1)^{k-j} \binom{k}{k-j} = \sum_{u=0}^{k-m} (-1)^u \binom{k}{u} = (-1)^{k-m} \binom{k-1}{k-m}$

(for the last equality see [3, p.82, identity 168]]).

In proving items 2 and 3, we use a well-known identity (see [3, p.65-66])

$$\sum_{j=0}^{\lfloor k/2 \rfloor} \binom{k}{2j} = \sum_{j=0}^{\lfloor (k-1)/2 \rfloor} \binom{k}{2j+1} = 2^{k-1}.$$

For 2, we have

$$c_k = \sum_{j=1}^{\lfloor k/2 \rfloor} (-1)^{k-2j} \binom{k}{2j} = (-1)^k \sum_{j=1}^{\lfloor k/2 \rfloor} \binom{k}{2j} = (-1)^k \left( \sum_{j=0}^{\lfloor k/2 \rfloor} \binom{k}{2j} - 1 \right) = (-1)^k (2^{k-1} - 1)$$

For 3, we have

$$c_k = \sum_{j=0}^{\lfloor (k-1)/2 \rfloor} (-1)^{k-(2j+1)} \binom{k}{2j+1} = (-1)^{k-1} \sum_{j=0}^{\lfloor (k-1)/2 \rfloor} \binom{k}{2j+1} = (-1)^{k-1} 2^{k-1}. \blacksquare$$